\documentclass[journal,a4paper]{IEEEtran}

\usepackage{cite}
\usepackage{epsfig}
\usepackage{amsmath,amssymb,fancybox}

\newcommand{\ld}{\lambda}

\newcommand{\rset}{\mathbb{R}}

\newcommand{\qed}{\hfill \qedbox \\[1ex]}
\newcommand{\xb}{\textbf{x}}

\newcommand{\bmtrx}{\left[\begin{array}}
\newcommand{\emtrx}{\end{array}\right]}

\newcommand{\qedbox}{
\begin{minipage}[b]{3mm}
{\unitlength1pt
\begin{picture}(7,7)
\thicklines  \put(-0.2,0){\line(1,0){6.7}}
\put(6,0){\line(0,1){6.2}}
\thinlines \put(0,0){\line(0,1){6}}
\put(0,6){\line(1,0){6}}
\end{picture}}
\end{minipage}}

\newtheorem{theorem}{\textbf{Theorem}}[section]
\newtheorem{assumption}[theorem]{\textbf{Assumption}}
\newtheorem{definition}[theorem]{\textbf{Definition}}
\newtheorem{lemma}[theorem]{\textbf{Lemma}}

\newtheorem{remark}[theorem]{\textbf{Remark}}

\newtheorem{algorithm}[theorem]{\textbf{Algorithm}}

\begin{document}

\title{ Application of a smoothing technique to decomposition
 in  convex optimization}

\author{Ion Necoara and Johan A.K. Suykens \thanks{The authors
are with the Katholieke Universiteit Leuven, Department of
Electrical Engineering,  ESAT, Kasteelpark Arenberg 10, B--3001
Leuven (Heverlee), Belgium,  tel/fax:+321632 0362/1970
(email:{\{ion.necoara,johan.suykens\}@esat.kuleuven.be}).} }

\markboth{IEEE Transactions on Automatic Control} {Necoara et al.:
Application of a smoothing technique to decomposition
 in  convex optimization}

\maketitle

\begin{abstract}
Dual decomposition is a powerful technique for deriving
decomposition schemes for convex optimization problems with
separable structure. Although the Augmented Lagrangian is
computationally more stable than the ordinary Lagrangian, the
\textit{prox-term} destroys the separability of the given problem.
In this paper we use another approach to obtain a smooth Lagrangian,
based on a smoothing technique developed  by Nesterov, which
preserves separability of the problem. With this approach we derive
a new decomposition method, called   \textit{proximal center
algorithm}, which from the viewpoint of efficiency estimates
improves the bounds on the number of iterations of the classical
dual gradient scheme by an order of magnitude. This can be achieved
with the new decomposition algorithm since  the resulting dual
function has good smoothness properties and since we  make use of
the particular structure of the given problem.
\end{abstract}





\section{Introduction}
\label{intro}

There has been  considerable recent interest in
parallel and distributed computation methods for solving large-scale
optimization problems  (e.g. \cite{BerTsi:89}). For separable convex
 problems, i.e. separable objective function but with coupling
constraints (this type of problems arise in many fields of
engineering: e.g. networks \cite{XiaBoy:04,PalChi:07}, distributed
model predictive control (MPC) \cite{VenRaw:07,NecDoa:cdc08},
stochastic programming \cite{RocWet:91},  etc), many researchers
have proposed dual decomposition algorithms such as the dual
subgradient method \cite{Uza:58,BerTsi:89},  alternating direction
method \cite{Tse:91,KonLeo:96,BerTsi:89,EckBer:92}, proximal method
of multipliers \cite{CheTeb:94}, partial inverse method
\cite{Spi:85,MahOua:95}, etc. In general, these methods are based on
alternating minimization in a Gauss-Seidel fashion of an (Augmented)
Lagrangian followed by a steepest ascent update  for the
multipliers. However, the step-size parameter which has a very
strong influence on the convergence rate of these methods is  very
difficult to tune and also they do not provide any complexity
estimates for the general  case (linear convergence is obtained e.g.
under strong convexity assumptions). Moreover, these methods use the
steepest ascent update for the multipliers, while we know from
\cite{Nes:04} that this update is inferior with  one order of
magnitude compared to Nesterov's accelerated scheme.

In this paper we  propose a new decomposition method for separable
convex optimization problems that overcomes  the disadvantages
mentioned above. Based on a smoothing technique recently developed
by Nesterov in \cite{Nes:05}, we obtain a smooth Lagrangian that
preserves separability of the problem. Using  this smooth
Lagrangian, we derive a new dual decomposition method in which the
corresponding parameters  are selected optimally and thus
straightforward to tune. In contrast to the dual gradient update for
the multipliers used  by most of the decomposition methods  from the
literature, our method uses an optimal gradient-based  scheme (see
e.g. \cite{Nes:04,Nes:05}) for updating the multipliers. Therefore,
we derive for the new method an efficiency estimate for the general
case which improves with one order of magnitude the complexity of
the classical dual gradient method (i.e. the steepest ascent
update). Up to our knowledge these are the first
 efficiency estimate results of a dual decomposition method for separable
  \textit{non-strongly} convex  programs.  The new algorithm is suitable for
decomposition since  it is highly parallelizable and thus it can be
effectively implemented on parallel processors. This is a distinct
feature of our method compared to alternating direction methods
based on Gauss-Seidel iterations that obviously do not share this
advantage.

 This paper is organized as follows. Section
 \ref{prel} contains the problem formulation, followed by a brief
introduction of some of the existing dual decomposition methods and
the description  of an
 accelerated scheme for smooth minimization  developed by Nesterov
 in \cite{Nes:04,Nes:05}. The main results of the paper are presented in
  Section \ref{newds}, where we describe our new
 decomposition method and its main properties, in particular global convergence.
  We conclude with some
  applications  and preliminary computational results on
  some test problems in Section \ref{section_apl}.


\section{Preliminaries}
\label{prel}

\subsection{Decomposition methods for separable convex programs}

 An important application of convex duality theory is in
decomposition algorithms  for solving large-scale problems but with
special structure. One such example, that we also consider in this
paper, is the following separable convex program:
\begin{align}
\label{scp} f^* = \min_{x \in X,z \in Z} \big \{ \phi_1(x) + \phi_2
(z): \; Ax + Bz = b \big \},
\end{align}
where $\phi_1:\rset^m \to \rset$ and $\phi_2: \rset^p \to \rset$ are
continuous convex functions on $X$ and $Z$, respectively, $A$ is a
given $n \times m$ matrix, $B$ is a given $n \times p$ matrix, and
$b$ is a given  vector in $\rset^n$. In this paper we do
\textit{not} assume $\phi_1$  and/or $\phi_2$ to be strongly convex
or smooth. Moreover, we assume that $X \subseteq \rset^m$ and  $Z
\subseteq \rset^p$ are given compact convex
 sets. We also use different norms on $\rset^n, \rset^m$ and
 $\rset^p$, not necessarily the corresponding Euclidian norms. However, for
 simplicity in notation we do not use indices to specify the norms
 on $\rset^n, \rset^m$ and
 $\rset^p$, since from the context it will be clear in which
 Euclidian space we are  and  which norm we use.
\begin{remark}
This type of problems  \eqref{scp} arises e.g. in the context of
large-scale networks   consisting  of multiple agents  with
different objectives or in the area of distributed model predictive
control (see also Section \ref{section_apl}).

\noindent We can also have any number $M$ of agents with different
$\phi_i$'s, not necessarily  two agents. Moreover, the
 method developed in this paper can handle  both  coupling  equalities ($Ax+Bz \!\!= \!\!b$) and/or
inequalities ($Ax+Bz \leq b$).
  However, for simplicity of the exposition we restrict  ourselves to
  \eqref{scp}. \qed
\end{remark}

Let   $\langle \cdot,\cdot \rangle$ denote the scalar product on the
Euclidian space $\rset^n$. By forming the Lagrangian corresponding
to the linear constraints (with the Lagrange multipliers $\ld \in
\rset^n$), i.e.
 $ \mathcal L_0(x,z,\ld) = \phi_1(x)+ \phi_2(z) + \langle \ld,Ax+Bz-b \rangle$,
  and using the dual decomposition method,
  one arrives at
the following decomposition algorithm:
\begin{algorithm}(\!\!\cite{Uza:58,BerTsi:89})
\label{alg1}
 \textbf{for} $k \geq 0$ \textbf{do}
\begin{enumerate}
\item[\textbf{1.}] given $\ld^k$, minimize the Lagrangian  $(x^{k+1},z^{k+1}) =
\arg \min_{x \in X, z \in Z} \mathcal L_0(x,z,\ld^k)$, or
equivalently minimize in parallel  over $x$ and $z$:
\begin{align*}
x^{k+1} & =\arg \min_{x \in X} [ \phi_1(x) + \langle \ld^k,Ax
\rangle], \\
z^{k+1} & =\arg \min_{z \in Z} [ \phi_2(z) + \langle \ld^k,Bz
\rangle ]
\end{align*}
\item[\textbf{2.}] update the multipliers: \[\ld^{k+1} = \ld^k + c_k (A x^{k+1} + B z^{k+1} -
b),\]
\end{enumerate}
\end{algorithm}
where $c_k$ is a positive step-size.\\
The following assumption is valid throughout the paper:
\begin{assumption}
\label{ass1} The set of optimal Lagrange multipliers $\Lambda^*$ is
nonempty for problem \eqref{scp}. \qed
\end{assumption}

It is known  that Algorithm \ref{alg1} is convergent under
Assumption \ref{ass1} and  the assumption that both $\phi_1$ and
$\phi_2$ are strongly convex functions (the latter guarantees that
the minimizer $(x^{k+1}, z^{k+1})$ is unique). In fact, under the
assumption of strong convexity,  the dual function
 \[f_0(\ld) = \min_{x \in X, z \in Z} \phi_1(x) + \phi_2(z) + \langle \ld,Ax+
 Bz-b \rangle\]
  is differentiable \cite{Roc:70,BerTsi:89}, and thus Algorithm \ref{alg1} can be seen  as the
 gradient method with step-size
 $c_k$ for maximizing the dual function.

 However, for many interesting problems, especially arising from transformations that
 leads to decomposition (see Section \ref{section_apl}),    the functions $\phi_1$ and
 $\phi_2$ are not strongly convex.  There are some methods (alternating direction
 method \cite{Tse:91,KonLeo:96,BerTsi:89}, proximal point method
  \cite{CheTeb:94}, partial inverse
 method \cite{Spi:85}) that
 overcome this difficulty based on e.g. alternating  minimization in a Gauss-Seidel fashion of the
  \textit{Augmented Lagrangian}, followed
by a steepest ascent update of the multipliers. A computational
drawback of these schemes is that the \textit{prox-term}
$\frac{c}{2}\|Ax + Bz -b \|^2$, using the Euclidian norm framework,
present in the Augmented Lagrangian is not separable in $x$ and $z$.
Another disadvantage is that they cannot deal with coupling
\textit{inequalities} in general. Moreover, these schemes   were
shown to be very sensitive to the value of the parameter $c$, with
difficulties in practice to obtain the best convergence rate. Some
heuristics for choosing $c$ can be found in the literature
\cite{KonLeo:96,Tse:91,CheTeb:94}. But, these heuristics have not
been formally analyzed from the viewpoint of efficiency estimates
for the general non-smooth case (linear convergence results were
obtained e.g. for strongly convex functions).
Note that alternating direction method variants which allow
for inexact minimization were proposed in
\cite{EckBer:92,CheTeb:94}. A closely related method is the partial
inverse of a monotone operator developed in \cite{Spi:85,MahOua:95}.


\subsection{An accelerated scheme  for smooth convex maximization}
\label{accel}

In this section we briefly describe  an accelerated scheme that also
uses only first-order information for smooth convex functions
developed by Nesterov  in \cite{Nes:04,Nes:05}.  Let $f$ be a
concave and differentiable function on a closed convex set $\hat Q
\subseteq \rset^n$. We further assume that the gradient of this
function is Lipschitz continuous:\\
 \hspace*{1cm}$\|\nabla f(x) -\nabla f(y) \|_* \leq L \| x - y \|  \;\;\; \forall x,y \in \hat Q$, \\
 where $\|s\|_* = \max_{\|x\| \leq 1} \langle s,x \rangle$ is the
corresponding dual norm of the norm used on $\rset^n$
\cite{Roc:70,BoyVan:04}.
\begin{definition} \cite{Nes:05}
\label{def_proxf}
 We define a \textit{prox-function} $d$ of the set
$\hat Q$ as a function with the following properties: \\
(i) $d$ is continuous, strongly convex on $\hat Q$ with convexity
parameter $\sigma$, \\
(ii) $u^0$ is the center of the set $\hat Q$, i.e. $u^0= \arg
\min_{x \in \hat Q} d(x)$ such that $d(u^0) = 0$.
\end{definition}
The goal is to find an approximate solution to the smooth convex
problem  $ x^* = \arg \max _{x \in \hat Q} f(x).$  In Nesterov's
scheme three sequences of points from $\hat Q$ are updated
recursively: $\{u^k\}_{k \geq 0}, \{x^k\}_{k \geq 0}$, and
$\{v^k\}_{k \geq 0}$. The algorithm can be described as follows:
\begin{algorithm} (\!\!\cite{Nes:05})
\label{alg4}
 \textbf{for} $k \geq 0$ \textbf{do}\\
\textbf{1.}  compute $f(u^k)$ and $\nabla f(u^k)$ \\
\textbf{2.} find
$x^k = \arg \max_{x \in \{\bar x^k, x^{k-1}, u^k \}}
f(x)$ where \\
$\bar x^k = \arg \max_{x \in \hat Q} [f(u^k) + \langle \nabla
f(u^k),x-u^k \rangle - \frac{L}{2} \| x-u^k \|^2]$ \\
\textbf{3.} find $v^k = \arg \max_{x \in \hat Q} \big \{-
\frac{L}{\sigma} d(x) + \sum_{l=0}^k \frac{l+1}{2} [f(u^l) + \langle
\nabla
f(u^l),x-u^l \rangle] \big \}$\\
\textbf{4.}  set $u^{k+1} = \frac{k+1}{k+3} x^k + \frac {2}{k+3}
v^k$.
\end{algorithm}
The derivation of Algorithm \ref{alg4} is based on the notion of
\textit{estimate sequence}. The main property of the estimate
sequence corresponding to Algorithm \ref{alg4} is the following
relation \cite{Nes:05}:
\begin{align}
\label{in_es} \frac{(k+1)(k+2)}{4} & f(x^k) \geq \max_{x \in \hat Q}
\big \{ - \frac{L}{\sigma} d(x) + \nonumber \\
&\sum_{l=0}^k \frac{l+1}{2}
[f(u^l) + \langle \nabla f(u^l),x-u^l \rangle] \big \}.
\end{align}
The convergence properties of Algorithm \ref{alg4} are summarized in
the following theorem:
\begin{theorem}\cite{Nes:05}
\label{thNes:05} Let   sequence $\{x^k\}_{k \geq 0}$  be generated
by Algorithm \ref{alg4}. Then,  $\{f(x^k)\}_{k \geq 0}$ is
nondecreasing  and we have the following efficiency estimate:
\[ f(x^*) - f(x^k) \leq  \frac{4 L d(x^*)}{\sigma (k+1)(k+2)}.\]
\end{theorem}
Theorem \ref{thNes:05} tells us that from the viewpoint of
efficiency estimates Nesterov's method applied to  maximization of a
concave function with Lipschitz continuous gradient has the order
$\mathcal O \big (\sqrt{\frac{L} {\epsilon}} \big )$. Therefore, the
efficiency of the method is higher by an order of magnitude  than
the corresponding pure gradient method (steepest ascent update with
complexity $ \mathcal O \big (\frac{1} {\epsilon} \big )$) for the
same  smooth problem (see \cite{Nes:04}). Note that we can define
directly $x^k = \bar x^k$ in step 2. The conclusions of Theorem
\ref{thNes:05} remain the same except that the sequence
$\{f(x^k)\}_{k \geq 0}$ is not necessarily monotone.


\section{A new decomposition method based on smoothing the Lagrangian}
\label{newds}

  In this section we propose a new method to smooth the Lagrangian
of \eqref{scp}, inspired from \cite{Nes:05}. This smoothing
technique preserves the separability of the problem   and moreover
the corresponding parameters are easy to tune.  Since separability
is preserved under this smoothing
 technique, we derive a new dual decomposition method in which  the multipliers are updated
 according to  Algorithm \ref{alg4}.  Moreover, we  obtain  efficiency estimates for the new
 method for the general  case and also global convergence. Note that with our method we can treat both  coupling
 equalities $Ax+Bz = b$ and/or inequalities $Ax+Bz \leq b$ (see also Remark \ref{cin}).

\subsection{Smoothing the Lagrangian}
\label{subs_sl}
 Let $d_X$ and $d_Z$  be two prox-functions for the compact convex sets $X$ and
 $Z$, with convexity parameter $\sigma_X$ and
 $\sigma_Z$, respectively.
 Denote  $x^0 =\arg \min_{x \in X} d_X(x), \;\; z^0 = \arg \min_{z \in Z} d_Z(z)$.
Since $X$ and $Z$ are compact and $d_X$ and $d_Z$ are continuous, we
can choose finite and positive constants
\[D_X \geq \max_{x \in X}
d_X(x), \; D_Z \geq \max_{z \in Z} d_Z(z). \]
 We also introduce the
following notation $\|A\| = \max_{\|\ld\|= 1,\|x\| = 1} \langle
\ld,A x \rangle$. Since the linear operator  $A$ is defined as $A:
\rset^m \to \rset^n_*$, where $\rset^n_*$ is the dual of $\rset^n$
(in fact $\rset^n_*=\rset^n$), we have
\[ \|A\| = \max_{\|x\| = 1}  \| A x\|_*   \; \text{and} \; \| A x\|_* \leq \|A \| \|x\| \; \forall x.  \]
Similarly for   $B$. Let us introduce the following family of
functions:
\begin{align}
\label{fc} f_c(\ld) = \min_{x \in X, z \in Z} [\phi_1(x)+\phi_2(z) +
& \langle \ld, Ax+Bz-b \rangle + \nonumber \\
& c \big ( d_X(x) + d_Z(z) \big )],
\end{align}
where $c$ is a positive \textit{smoothness parameter} that will be
defined later in this section. Note that we could also choose
different parameters $c_1$ and $c_2$ for each prox-term. The
generalization is straightforward.   It is clear that the objective
function in \eqref{fc} is separable in $x$ and $z$, i.e.
\begin{align}
\label{fc1} f_c(\ld) =  - \langle \ld,b \rangle + & \min_{x \in X}
[\phi_1(x) + \langle \ld, Ax \rangle + c \ d_X(x)] + \nonumber \\
& \min_{ z \in Z} [ \phi_2(z) + \langle \ld, Bz \rangle + c \
d_Z(z)].
\end{align}
Denote by $x(\ld)$ and $z(\ld)$ the optimal solution  of the
minimization problem in $x$ and $z$, respectively.  Function $f_c$
has the following  smoothness properties:
\begin{theorem}
\label{th_fc} The function $f_c$ is concave and  continuously
differentiable at any $\ld \in \rset^n$. Moreover, its gradient
$\nabla f_c(\ld) = A x(\ld) + B z(\ld) -b$ is Lipschitz continuous
with Lipschitz constant $L_c=\frac{\|A\|^2}{c \sigma_X} +
\frac{\|B\|^2}{c \sigma_Z}$. The following inequalities  hold:
\begin{equation}
\label{f0fc}
 f_c(\ld) \geq f_0(\ld) \geq f_c(\ld) - c(D_X+D_Z) \;\; \forall \ld \in
\rset^n.
\end{equation}
\end{theorem}
\IEEEproof Since the functions $d_X$ and $d_Z$ are strongly convex,
it follows that the optimal solution $(x(\ld),z(\ld))$ of \eqref{fc}
or
 \eqref{fc1} is unique for any $\ld$ and thus the function $f_c$ is well defined
at any~$\ld$. Concavity and continuous differentiability of $f_c$
follows from standard duality theory \cite{Roc:70,BerTsi:89}. It
remains to show that its gradient  $\nabla f_c(\ld) = A x(\ld) + B
z(\ld) -b$ is Lipschitz continuous.   For simplicity of notation we
assume that all  the functions involved in the minimization problem
\eqref{fc} are  differentiable. Let $\ld$ and $\eta$ be two Lagrange
multipliers. Using first-order optimality conditions for the
minimization problem in $x$  we obtain:
\begin{align*}
&  \langle  \nabla \phi_1(x(\ld)) + A^T \ld + c \nabla d_X(x(\ld),
x(\eta)
-x(\ld) \rangle \geq 0\\
&  \langle \nabla \phi_1(x(\eta)) + A^T \eta + c \nabla
d_U(x(\eta)), x(\ld) -x(\eta) \rangle \geq 0.
\end{align*}
Adding these two inequalities and since $\phi_1$ is convex and $d_X$
is strongly convex, we obtain
\begin{align*}
& \langle  A^T(\eta - \ld) , x(\ld) -x(\eta) \rangle
 \geq \\
&  \langle \nabla \phi_1(x(\ld)) - \nabla \phi_1(x(\eta)), x(\ld) -
x(\eta)  \rangle
 + \\
 & c \langle \nabla d_X(x(\ld)) - \nabla d_X(x(\eta)), x(\ld) - x(\eta)  \rangle
 \geq \\
&   c \sigma_X \|x(\ld) - x(\eta) \|^2.
\end{align*}
From last relation and  Cauchy-Schwartz  inequality we have:
\begin{align*}
 \| A x(\ld) -A x(\eta)  \|_*^2 & \leq  \|A\|^2 \| x(\ld) - x(\eta)
\|^2  \\
 & \leq  \frac{\|A\|^2}{c \sigma_X} \langle A^T(\eta - \ld),
x(\ld) -x(\eta) \rangle \\
& \leq   \frac{\|A\|^2}{c \sigma_X} \| \ld - \eta \| \|A x(\ld) - A
x(\eta)\|_*,
\end{align*}
and thus $\| A x(\ld) -A x(\eta)  \|_* \leq \frac{\|A\|^2}{c
\sigma_X} \| \ld - \eta \|$.  Similarly, for the minimization
problem in $z$ we obtain: $ \|B z(\ld) -B z(\eta) \|_*  \leq
\frac{\|B\|^2}{c \sigma_Z} \| \ld - \eta \|$.  In conclusion, the
gradient of $f_c$ satisfies
\[  \| \nabla f_c(\ld) -\nabla f_c(\eta)  \|_*  \leq \big ( \frac{\|A\|^2}{c
\sigma_X} + \frac{\|B\|^2}{c \sigma_Z} \big ) \| \ld - \eta\|. \]
Furthermore, the first inequality in \eqref{f0fc} is a consequence
of the fact that $d_X(x) \geq 0$ for all $x$, and $d_Z(z) \geq 0$
for all $z$. The second inequality in \eqref{f0fc} follows from:
$f_c (\ld) \leq \min_{x \in X, z \in Z} [ \phi_1(x)+\phi_2(z) +
\langle \ld, Ax+Bz-b \rangle] + c  \max_{x \in X, z \in Z}  [ d_X(x)
+ d_Z(z)]$. \qed

\subsection{A proximal center--based decomposition method}

In this section we derive a new dual decomposition method  based on
the smoothing technique  described in Section \ref{subs_sl}. The new
algorithm,  called here the \textit{proximal center algorithm}, has
the nice feature that
 the \textit{coordination} between the agents involves the maximization of a
 smooth convex objective function (i.e.  with Lipschitz continuous
 gradient). Moreover, the \textit{resource allocation} stage consists in
 solving in \textit{parallel} two independent minimization problems
 with strongly convex objectives. The new method belongs to the class of two-level algorithms
 \cite{Coh:78} and is particularly
 suitable for  separable  convex problems where the minimizations
 over $x$ and $z$ in \eqref{fc1} are easily carried out.

\noindent We apply the accelerated method described in Algorithm
\ref{alg4} to the concave function $f_c$ with Lipschitz continuous
gradient:
\begin{equation}
\label{min_fc} \max_{\ld \in Q} f_c(\ld),
\end{equation}
where  $Q$ is a given closed convex set in $\rset^n$ that contains
at least one optimal multiplier $ \ld^* \in \Lambda^*$. Notice that
$Q \subseteq \rset^n$ for  linear equalities (i.e. $Ax + Bz -b=0$),
$Q \subseteq \rset^n_+$, where $\rset_+$ denotes the set of
nonnegative real numbers, for linear inequalities (i.e. $A x + Bz -b
\leq 0$), or $Q \subseteq (\rset^{n_1} \times \rset^{n_2}_+)$, where
$n_1+n_2 =n$, when both, equalities and inequalities are present.
Note that according to Algorithm \ref{alg4} we also need to choose a
prox-function $d_Q$ for the set $Q$ with the convexity parameter
$\sigma_Q$  and center $u^0$. The \textit{proximal center algorithm}
can be described as follows:
\begin{algorithm}
\label{alg5}
\textbf{for} $k \geq 0$ \textbf{do}
\begin{enumerate}
\item[\textbf{1.}] given $u^k$ compute in \textit{parallel}
\[x^{k+1} = \arg \min_{x \in X} [\phi_1(x) + \langle u^k, Ax \rangle + c \
d_X(x)] \] \[ z^{k+1} = \arg \min_{ z \in Z} [ \phi_2(z) + \langle
u^k, Bz \rangle + c \ d_Z(z)]. \]
\item[\textbf{2.}] compute
\[f_c(u^k) \!=\! \mathcal L_0(x^{k\!+\!1}\!,z^{k\!+\!1}\!,u^{k})
+ c \big ( d_X(x^{k\!+\!1}) + d_Z(z^{k\!+\!1}) \big ),\]
\[\nabla f_c(u^k)= A x^{k+1} \!+\! B z^{k+1} - b \]
\item[\textbf{3.}] find $\ld^k = \arg  \max_{\ld \in \{\bar \ld^k,
\ld^{k-1}, u^k \}} f_c(\ld)$ where
\[ \bar \ld^k \!=\! \arg \max_{\ld \in
Q} f_c(u^k) \!+\! \langle \nabla f_c(u^k),\ld-u^k \rangle \!-
\frac{L_c}{2} \| \ld-u^k \|^2\]
\item[\textbf{4.}] find
\begin{align*}
v^k = \arg & \max_{\ld \in  Q} \big \{  - \frac{L_c}{\sigma_Q}
d_Q(\ld)
+ \\
& \sum_{l=0}^k \frac{l+1}{2} [f_c(u^l) + \langle \nabla
f_c(u^l),\ld-u^l \rangle] \big \}
\end{align*}
\item[\textbf{5.}] set $u^{k+1} = \frac{k+1}{k+3} \ld^k + \frac {2}{k+3} v^k$.
\end{enumerate}
\end{algorithm}
The proximal center algorithm is  suitable for decomposition since
it is highly parallelizable: the  agents can solve in parallel their
corresponding minimization problems. We now derive a lower bound for
the value   of the objective function which will be used frequently
in the sequel:
\begin{lemma}
\label{lowerpo} For any $\lambda^* \in \Lambda^*$ and $\hat x \in X,
\hat z \in Z$, the following lower bound on primal  gap holds:
\[ [\phi_1(\hat x) + \phi_2(\hat z)] - f^* \geq  - \|\ld^*\|_* \|A \hat x + B \hat z -
b\|_*. \]
\end{lemma}
\IEEEproof From our assumptions we have that
\begin{align*}
f^* = & \min_{x \in X, z \in Z} [\phi_1(x) + \phi_2(z) -  \langle A
x + B z - b,\ld^* \rangle]\\
&  \leq \phi_1(\hat x) + \phi_2(\hat z)
- \langle A \hat x + B \hat z - b,\ld^* \rangle,
\end{align*}
and then using the Cauchy-Schwarz  inequality, the result
follows.\qed
 The previous lemma shows that if $\| A \hat x + B \hat z - b
\|_* \leq \epsilon_c$, then the primal gap is bounded: for all $\hat
\ld \in Q$
\begin{align}
\label{lowerpo1}
 -\epsilon_c \|\ld^*\|_* \leq  \phi_1(\hat
x) + \phi_2(\hat z)]&  - f^*  \leq \nonumber \\
& \phi_1(\hat x) + \phi_2(\hat z) - f_0 (\hat \ld).
\end{align}
Therefore, if we are able to derive an upper bound  $\epsilon$  for
the duality gap and $\epsilon_c$ for the coupling constraints for
some given $\hat \ld$ and $\hat x \in X, \hat z \in Z$, then we
conclude  that $(\hat x, \hat z)$ is an $(\epsilon,
\epsilon_c)$-solution for problem \eqref{scp} (since in this case
$-\epsilon_c \|\ld^*\|_* \leq \phi_1(\hat x) + \phi_2(\hat z) - f^*
\leq \epsilon$ for all $\ld^* \in \Lambda^*$). The next theorem
derives an upper bound on the duality gap for our method.
\begin{theorem}
\label{th_conv}
 Assume that there exists a closed convex set $Q$
that contains a $\ld^* \in \Lambda^*$. Then, after $k$ iterations we
obtain an approximate solution to the problem \eqref{scp} \[ (\hat
x,\hat z) = \sum_{l=0}^k \frac{2(l+1)}{(k+1)(k+2)}
(x^{l+1},z^{l+1})\] and $\hat \ld=\ld^k$ which satisfy the following
duality gap:
\begin{align}
\label{ubdg}
& [\phi_1(\hat x) + \phi_2(\hat z)] - f_0 (\hat \ld)
\leq c(D_X+D_Z) - \nonumber \\
& \max_{\ld \in Q} \big [ -\frac{4 L_c}{\sigma_Q(k+1)^2} d_Q(\ld) +
 \langle A \hat x + B \hat z - b,\ld \rangle \big ].
\end{align}
\end{theorem}
\IEEEproof For an arbitrary $c$, we have from the inequality
\eqref{in_es} that after $k$ iterations the Lagrange multiplier
$\hat \ld$ satisfies  the following relation:
\begin{align*}
 & \frac{(k+1)(k+2)}{4} f_c(\hat \ld) \geq  \max_{\ld
\in Q} \big \{ - \frac{L_c}{\sigma_Q} d_Q(\ld) + \nonumber \\
 &\sum_{l=0}^k \frac{l+1}{2} [f_c(u^l) +
\langle \nabla f_c(u^l),\ld-u^l \rangle] \big \}.
\end{align*}
In view of  the previous inequality  we have:
\begin{align*}
f_c(\hat \ld) \geq  &  \max_{\ld \in Q}  \big \{ - \frac{4 L_c}{\sigma_Q (k+1)^2} d_Q(\ld)  + \\
 &  \sum_{l=0}^k \frac{2(l+1)}{(k+1)(k+2)} [f_c(u^l) +
\langle \nabla f_c(u^l),\ld-u^l \rangle] \big \}.
\end{align*}
 Now, we replace $f_c(u^l)$  and $\nabla f_c(u^l)$ with the
 expressions given in step 2 of Algorithm \ref{alg5} we obtain: for all
$\ld \in Q$
\begin{align*}
& \sum_{l=0}^k \frac{2(l+1)}{(k+1)(k+2)} [f_c(u^l) + \langle \nabla
f_c(u^l),\ld-u^l \rangle] \geq \\
& \sum_{l=0}^k \frac{2(l+1)}{(k+1)(k+2)} [ \langle A x^{l+1} + B
z^{l+1} -
b,\ld \rangle+ \\
& \phi_1(x^{l+1}) + \phi_2(z^{l+1}) ] \geq  \langle A \hat x + B
\hat z - b,\ld \rangle + \phi_1(\hat x) + \phi_2(\hat z).
\end{align*}
The first inequality follows from the fact that the prox-functions
$d_X, d_Z \geq 0$ and the last inequality follows from convexity of
the functions $\phi_1$ and $\phi_2$. Using the last relation and
\eqref{f0fc} we derive the  bound \eqref{ubdg} on the duality gap.
\qed
 We now show how to construct the set $Q$
and how to   choose optimally the smoothness parameter $c$.  In the
next two sections we discuss two cases depending on the choices for
$Q$ and $d_Q$.

\subsection{Efficiency estimates for compact $Q$}

Let $D_Q$ be a positive constant satisfying
\begin{align}
\label{max_dq}
 \max_{\ld \in Q} d_Q(\ld) \leq D_Q.
\end{align}
Let us note that we can choose  $D_Q$ finite whenever $Q$ is
compact. In this section we specialize the result of Theorem
\ref{th_conv} for the case when $Q$ has the following form:  \[Q=\{
\ld \in \rset^n: \| \ld \| \leq R \}.\]
\begin{theorem}
\label{th_b} Assume that $\Lambda^*$ is  bounded. Then, the sequence
$\{\ld^k\}_{k \geq 0}$ generated by Algorithm~\ref{alg5} is also
bounded.
\end{theorem}
\IEEEproof  Note that $\Lambda^*=\{ \ld : f_0(\ld) \geq f^* \}$. Let
us introduce the sets $ \Lambda^0 = \{ \ld : f_0(\ld) \geq f^* - c
(D_X+D_Z) \}$  and
 $\Lambda^c=\{\ld : f_c(\ld) \geq f^* \}$. From the inequalities in \eqref{f0fc} it follows
 immediately that
$\Lambda^* \subseteq \Lambda^0$ and also $\Lambda^* \subseteq
\Lambda^c$. Therefore, the sets $\Lambda^c$ and $\Lambda^0$ are
nonempty. Since $\Lambda^*$ is bounded, from Corollary 8.7.1 in
\cite{Roc:70} it follows that the set $\Lambda^0$ is also bounded.
We can also show that  $\Lambda^c$ is a bounded set. Indeed let $\ld
\in \Lambda^c$, then  using once more the second inequality in
\eqref{f0fc}  we obtain: $f_0(\ld) + c (D_X + D_Z) \geq f_c(\ld)
\geq f^*$, i.e. $\ld \in \Lambda^0$. In conclusion,  $\Lambda^c
\subseteq \Lambda^0$ and thus $\Lambda^c$ is also bounded.

Let us now show that the  sequence $\{\ld^k\}_{k \geq 0}$ is
bounded. From Theorem \ref{thNes:05} it follows that the sequence
$\{f_c(\ld^k)\}_{k \geq 0}$ is nondecreasing and thus $\{ \ld^k: k
\geq 0 \}  \subseteq \{ \ld: f_c(\ld) \geq f_c(\ld^0) \}$. But since
$\Lambda^c$ is bounded, using once again Corollary 8.7.1 in
\cite{Roc:70} it follows that the set $\{ \ld: f_c(\ld) \geq
f_c(\ld^0) \}$ is bounded. In conclusion, the sequence $\{\ld^k\}_{k
\geq 0}$ is bounded. \qed Since  Assumption \ref{ass1} holds, then
$\Lambda^*$ is nonempty. Conditions under which $\Lambda^*$ is
bounded can be found in \cite{Roc:70} (e.g. when the matrix $[A \;
B]$ has full  rank). Under the assumptions of Theorem \ref{th_b}, it
follows that there exists $R>0$ sufficiently large such that the set
$Q = \{\ld \in \rset^n: \|\ld \| \leq R\}$ contains\footnote{In a
practical algorithm $R$
 is increased adaptively: if $\ld^k$ approaches the boundary of $Q$ we  take
 $R_{k+1} = \alpha R_k$ for some $\alpha>1$. An upper bound on $\|\ld^*\|_*$
 can also  be  estimated using $R$.} a  $\ld^*
\in \Lambda^*$, and thus we can assume $D_Q$ to be finite. Notice
that similar arguments were used in order to prove convergence of
two-level algorithms for convex  problems in \cite{Coh:78}.

The next theorem shows how to choose optimally the smoothness
parameter $c$ and provides  the complexity estimates of  our method
for the case when $Q$ is a ball.
\begin{theorem}
\label{th_conv1}
 Assume that there exists $R>0$ such that the  set $Q =\{\ld \in \rset^n: \|\ld
\| \leq R\}$  contains a $\ld^* \in \Lambda^*$. Taking $c =
\frac{2}{k+1} \sqrt{
 \frac{D_Q}{D_X+D_Z} \big ( \frac{\|A\|^2}{\sigma_Q \sigma_X} +
 \frac{\|B\|^2}{\sigma_Q \sigma_Z} \big ) }$,
  then after $k$ iterations we obtain an approximate solution to
the problem \eqref{scp} $(\hat x,\hat z) = \sum_{l=0}^k
\frac{2(l+1)}{(k+1)(k+2)} (x^{l+1},z^{l+1})$ and $\hat \ld=\ld^k$
which satisfy the following bounds on the  duality gap and
constraints:
\begin{align*}
& [\phi_1(\hat x) + \phi_2(\hat z)] - f_0 (\hat \ld)
\leq \nonumber \\
& \frac{4} {k+1} \sqrt{D_Q(D_X+D_Z) \big ( \frac{\|A\|^2}{\sigma_Q
\sigma_X} +
 \frac{\|B\|^2}{\sigma_Q \sigma_Z} \big ),}
\end{align*}
\begin{align*}
&  \| A \hat x + B \hat z - b \|_* \leq  \\
 &  \frac{4}
{(R-\|\ld^*\|_*)(k+1)} \sqrt{D_Q(D_X+D_Z) \big (
\frac{\|A\|^2}{\sigma_Q \sigma_X} + \frac{\|B\|^2}{\sigma_Q
\sigma_Z} \big )}.
\end{align*}
\end{theorem}
\IEEEproof Using \eqref{max_dq} and the  form of $Q$ we  obtain that
\begin{align*}
& \max_{\ld \in Q}  -\frac{4 L_c}{\sigma_Q(k+1)^2} d_Q(\ld)  +
\langle A \hat x + B \hat z - b,\ld \rangle \geq \\
&  - \frac{4 L_c
D_Q}{\sigma_Q (k+1)^2}  +    \max_{\| \ld \| \leq R } \langle A \hat
x + B \hat z - b,\ld
\rangle = \\
& - \frac{4 L_c D_Q}{\sigma_Q (k+1)^2}  +  R  \| A \hat x + B \hat z
- b\|_*.
\end{align*}
In view of the previous relation and  Theorem \ref{th_conv} we
obtain the  following duality gap:
\begin{align*}
 & [\phi_1(\hat x) + \phi_2(\hat z)] - f_0 (\hat \ld) \leq
c(D_X+D_Z) +\frac{4 L_c D_Q}{\sigma_Q (k+1)^2}  - \\
& R  \| A \hat x + B \hat z - b\|_* \leq  c(D_X+D_Z) +\frac{4 L_c
D_Q}{\sigma_Q (k+1)^2}.
\end{align*}
Minimizing the right-hand side of this inequality over $c$ we get
the  above expressions for $c$ and for the upper bound on the
duality gap. Moreover, for the constraints using Lemma \ref{lowerpo}
and  inequality \eqref{lowerpo1} we have that
\[ (R-\|\ld^*\|_*)  \| A \hat x + B \hat z - b\|_*    \leq
 c(D_X+D_Z) +\frac{4 L_c D_Q}{\sigma_Q (k+1)^2} \]
and replacing  $c$ derived above we also get the upper bound on the
constraints violation. \qed
From Theorem \ref{th_conv1} and
inequality \eqref{lowerpo1} we obtain that the complexity for
finding an $(\epsilon,\epsilon_c)$-approximation of the optimal
value function $f^*$, when the set $Q$ is a ball, is $k+1=4
\sqrt{D_Q(D_X+D_Z) \big ( \frac{\|A\|^2}{\sigma_Q \sigma_X} +
 \frac{\|B\|^2}{\sigma_Q \sigma_Z} \big )}\; \frac{1} {\epsilon}$,
 i.e.  the efficiency estimates of our scheme is of the order ${\cal
 O}(\frac{1}{\epsilon})$, better than most non-smooth optimization
 schemes such as the subgradient method that have an efficiency estimate of the order ${\cal
 O}(\frac{1}{\epsilon^2})$ (see e.g. \cite{Nes:04}).  Moreover the
 dependence of the parameters $c$ and $L_c$  on $\epsilon$ is as
 follows: $c= \frac{\epsilon}{2(D_X+D_Z)}$ and $L_c=\big (
 \frac{\|A\|^2}{\sigma_X} +  \frac{\|B\|^2}{\sigma_Z} \big ) \frac{D_X+D_Z}{2
\epsilon}$. Another advantage  of the proximal center method is that
we are free in the choice of the norms in the spaces  $\rset^n,
\rset^m$ and $\rset^p$, while most
 of the decomposition schemes are based on the Euclidian norm. Thus, we can
 choose the norms which make  the ratio $\frac{\|A\|^2}{\sigma_Q \sigma_X}$
 as small as possible.

\subsection{Efficiency estimates for the Euclidian norm}
\label{sec_en}

In this section we assume that $\rset^n$ is endowed with the
Euclidian norm.
\begin{theorem}
\label{th_conv2}
 Assume that $Q=\rset^n$ and $d_Q(\ld) = \frac{1}{2}
\|\ld\|^2$, with the Euclidian norm on $\rset^n$. Taking $c=
\frac{\epsilon}{D_X+D_Z}$ and\\ $k+1 = 2 \sqrt{ \big
(\frac{\|A\|^2}{\sigma_X} +
 \frac{\|B\|^2}{\sigma_Z} \big )(D_X+D_Z)} \; \frac{1}{\epsilon}$, then
 after $k$ iterations the duality gap is less than $\epsilon$ and
 the constraints satisfy  $\|A \hat x + B \hat z - b\|
\leq  \epsilon \big ( \|\ld^*\| + \sqrt{\|\ld^*\|^2+2} \big )$.
\end{theorem}
\IEEEproof Let us note that for these choices of $Q$ and $d_Q$ we
have $\sigma_Q=1$ and thus
\begin{align*}
 & \max_{\ld \in \rset^n}  -\frac{4
L_c}{\sigma_Q (k+1)^2} d_Q(\ld)  + \langle A \hat x + B \hat z -
b,\ld \rangle =\\
& \frac{(k+1)^2}{8 L_c} \|A \hat x + B \hat z - b\|^2.
\end{align*}
In this case we obtain the following  bound on the duality gap (see
Theorem \ref{th_conv}):
\begin{align*}
 & [\phi_1(\hat x) + \phi_2(\hat z)] - f_0 (\hat \ld) \leq \\
 & c(D_X+D_Z)
- \frac{(k+1)^2}{8 L_c} \|A \hat x + B \hat z - b\|^2 \leq
c(D_X+D_Z).
\end{align*}
It follows that taking $c = \frac{\epsilon}{D_X+D_Z}$, the duality
gap is less than $\epsilon$.  For the constraints using Lemma
\ref{lowerpo} and inequality \eqref{lowerpo1} we get that $\|A \hat
x + B \hat z - b\|$ satisfies the second order inequality in $y$:  $
\frac{(k+1)^2}{8 L_c} y^2  - \|\ld^*\| y -\epsilon \leq 0$.
Therefore, $\|A \hat x + B \hat z - b\|$ must be less than the
largest root of the corresponding second-order equation, i.e.
\[\|A \hat x + B \hat z - b\|
\leq \big ( \| \ld^*\| + \sqrt{\|\ld^*\|^2 +
\frac{\epsilon(k+1)^2}{2 L_c}} \big ) \frac{4 L_c}{(k+1)^2}. \] With
some straightforward computations we get that after $k$ iterations,
where $k$ defined as in the theorem, we
 also get  $\|A \hat x + B \hat z - b\| \leq
  \epsilon \big (\|\ld^*\| + \sqrt{\|\ld^*\|^2+2} \big )$. \qed

\begin{remark}
\label{cin} (i) When coupling inequalities $A x + B z - b \leq 0$
are present, then we choose $Q = \rset^n_+$. Using the same
reasoning as before we get that $\max_{\ld \geq 0} -\frac{4
L_c}{\sigma_Q (k+1)^2} d_Q(\ld)  + \langle A \hat x + B \hat z -
b,\ld \rangle = \frac{(k+1)^2}{8 L_c} \|[A \hat x + B \hat z -
b]^+\|^2$, where $[y]^+$ denotes the projection of $y \in \rset^n$
onto $\rset^n_+$. Taking for $c$ and $k$ the same values as in
Theorem \ref{th_conv2}, we can conclude that after $k$ iterations
the duality gap is less than $\epsilon$ and using a modified version
of Lemma \ref{lowerpo} (i.e. a generalized version of Cauchy-Schwarz
inequality: $-\langle y, \ld
 \rangle \geq - \|[y]^+\| \ \|\ld\|$ for any $y \in \rset^n, \ld \geq 0$)
 the constraints violation satisfy $\|[A \hat x + B \hat z -
b]^+\| \leq \epsilon \big (\|\ld^*\| + \sqrt{\|\ld^*\|^2+2} \big )
$.\\
 (ii) Note that our algorithm can also deal with more general
inequality constraints, e.g. sum of separable convex functions.
\end{remark}

Finally, let us
 mention  that our decomposition method described in Algorithm
 \ref{alg5} bears similarity with  proximal type methods
 \cite{BerTsi:89,EckBer:92,Tse:91,CheTeb:94,MahOua:95,KonLeo:96},
 but is different both in the computational steps and in the choice
 of the parameter $c$.  More precisely,  our method uses
  a fixed center in the prox-terms  that allows
 the inner  iterations at each $k$ to move freely, contrary to most
 proximal type methods that force the next iterates to be close to
 the previous ones. The main advantage of our scheme is that  it
 is fully automatic, the parameter $c$ is chosen unambiguously, which
 is crucial for justifying the convergence properties of  Algorithm
 \ref{alg5}.


\section{Applications}
\label{section_apl}

\subsection{Applications with separable structure}
\label{subsec_apl} In this section we briefly  discuss some  of the
applications to which our method can be applied.

First application that we will discuss  here is the control of
large-scale systems with interacting subsystem dynamics. A
distributed model predictive control (MPC) framework is appealing
in this context since this framework allows us to design local
 controllers that take care of the interactions
between different subsystems and physical constraints. We assume
that the overall system model can be decomposed  into $M$
appropriate subsystem models:
\begin{align*}
x^i(k+1) = \sum_{j \in {\mathcal N}(i)} A_{ij}  x^j(k) + B_{ij}
u^j(k) \;\;  \forall i=1 \cdots M,
\end{align*}
where ${\mathcal N}(i)$ denotes the set of subsystems that interact
with the $i$th subsystem, including itself. The control and state
sequence  must satisfy local constraints:  $x^i(k) \in \Omega_i$ and
$u^i(k) \in U_i$  for all $i$ and  $k \geq 0$, where the sets
$\Omega_i$ and $U_i$ are usually convex compact sets with the origin
in their interior. In general the control objective is to steer the
state of the system to origin or any other set point in a ``best''
way. Performance is expressed via a stage cost, which in this paper
we assume to have the following form \cite{VenRaw:07}: $
\sum_{i=1}^M \ell_i(x^i,u^i)$, where usually $\ell_i$ is a convex
quadratic function, not necessarily strictly convex. Let $N$ denote
the prediction horizon. In MPC we must solve at each step $k$, given
$x^i(k) = x^i$, an optimal control problem of the following form:
\begin{align}
\label{dmpc}
 & \min_{x_l^i,u_l^i} \big \{  \sum_{l=0}^{N-1} \sum_{i=1}^M
\ell_i(x^i_l,u^i_l):   \; x_l^i \in \Omega_i, \; u_l^i \in U_i \;
\forall \; l, i \big \}
\end{align}
where $x_0^i=x^i$ and $x_{l+1}^i =\sum_{j \in {\mathcal N}(i)}
A_{ij} x^j_l + B_{ij} u^j_l$. A similar formulation of distributed
MPC for coupled linear subsystems with decoupled costs was given in
\cite{VenRaw:07}, but without  state constraints. Let us introduce
$\xb_i =(x_0^i \cdots x_N^i \ u_0^i \cdots u_{N-1}^i), X_i =
\Omega_i^{N+1} \times U_i^N$ and the non-strictly convex quadratic
functions $\psi_i(\xb^i) = \sum_{l=0}^{N-1} \ell_i(x^i_l,u^i_l)$.
Then, the control problem \eqref{dmpc} can be recast as a separable
convex program:
\begin{align}
\label{edmpc}
 \min_{\xb_i \in X_i}  \big \{
\sum_{i=1}^M \psi_i(\xb_i): \sum_{i=1}^M C_i \xb_i -\gamma =0 \big
\},
\end{align}
where  the $C_i$'s and  $\gamma$ are defined appropriately.

Network optimization  furnishes another areas in which our Algorithm
\ref{alg5}  leads to a new method of solution. In this application
the convex problem has the form \cite{XiaBoy:04,PalChi:07}:
\begin{align}
\label{general_dcop}
 \min_{x_i \in X_i}  \big \{
\sum_{i=1}^M \psi_i(x_i): \sum_{i=1}^M C_i x_i -\gamma =0, \;
\sum_{i=1}^M D_i x_i - \beta \leq 0  \big  \},
\end{align}
 where $X_i$ are compact convex sets (in general balls) in the Euclidian
 space $\rset^m$, $\psi_i$'s are non-strictly convex functions and $M$  denotes the number of agents
  in the network.

In \cite{VenRaw:07} the  optimization problem \eqref{dmpc} (or
equivalently \eqref{edmpc}) was  solved    in a decentralized
fashion, iterating the Jacobi algorithm  \cite{BerTsi:89} $p_{\max}$
times. But, there is no theoretical guarantee of the Jacobi
algorithm about how good is the approximation of the optimum after
$p_{\max}$ iterations and moreover we need strictly convex functions
$\psi_i$ to prove asymptotic convergence to the optimum. However, if
we solve \eqref{dmpc} using  Algorithm \ref{alg5} (see
\cite{NecDoa:cdc08} for more details), we have a guaranteed upper
bound (see Theorem \ref{th_conv1} or \ref{th_conv2}) on the
approximation of the optimum after $p_{\max}$ iterations. In
\cite{XiaBoy:04,PalChi:07} the optimization problem
\eqref{general_dcop} is solved using the dual subgradient method
described in the Algorithm \ref{alg1}. Some preliminary simulation
tests from Section \ref{simulations} show that  Algorithm \ref{alg5}
is superior to Algorithm \ref{alg1}.

Let us  describe briefly  the main ingredients of Algorithm
\ref{alg5} for problem \eqref{general_dcop}.  Let $d_{X_i}$  be
properly chosen  prox-functions, according to the structure  of the
sets $X_i$'s and the norm\footnote{For example if $X_i = \{ x \in
\rset^m : \|x - x_0\|_2 \leq R \}$, where  $\| \cdot \|$ denotes
here the Euclidian norm, then it is natural to take $d_{X_i} (x) =
\frac{\|x - x_0\|_2^2}{2}$. If $X_i = \{x \in \rset^m : x \geq 0, \;
\sum_{i=1}^m x_{(i)} =1 \}$, then the norm $\|x\| = \sum_{i=1}^m
|x_{(i)}|$ and $d_{X_i} (x) = \ln m + \sum_{i=1}^m x_{(i)} \ln
x_{(i)}$ is more suitable (see \cite{Nes:05} for more details).}
used on $\rset^m$. Then, the smooth dual function $f_c$ has the
form:
\begin{align*}
f_c(\ld) =  \min_{x_i  \in X_i} & \sum_{i=1}^M \psi_i(x_i) + \langle
\ld_1,   \sum_{i=1}^M C_i x_i -\alpha \rangle + \\
& \langle \ld_2, \sum_{i=1}^M D_i x_i -\beta \rangle +
 c \sum_{i=1}^M  d_{X_i} (x_i).
\end{align*}
Moreover, $Q \subseteq \rset^{n_1} \times \rset^{n_2}_+$, where
$n_1$ and $n_2$ denote the number of equality and inequality
constraints. In this case all the minimization problems of Algorithm
\ref{alg5} are decomposable in $x_i$ and thus the agents can solve
the optimization problem \eqref{general_dcop} in a distributed
fashion.

\subsection{Computational results}
\label{simulations}

 We conclude this paper  with the results  of
computational experiments on a random set  of problems  of the form
\eqref{general_dcop}, where $\psi_i$'s are convex quadratic
functions (but not strictly convex) and $X_i$'s are balls in
$\rset^m$ defined with the Euclidian norm. Similarly, $\rset^n$
(corresponding to the Lagrange multipliers) is endowed with the
Euclidian norm (see Section \ref{sec_en}).

\begin{table}[!h]
\label{table_sim}
\begin{center}
   \begin{tabular}{| c | c || c | c || c| c| }
     \multicolumn{2}{c}{} &   \multicolumn{2}{c}{$M=2$}         & \multicolumn{2}{c}{$M=10$}\\
     \hline
       m    & $\epsilon$  &       PCM       &    DSM            &     PCM         &    DSM   \\   \hline
       50   &    0.01     &      202        &  5000(0.05)       &     811        &   5000(0.29)       \\   \hline
      200   &    0.01     &      625        &  5000(0.19)       &     2148        &   5000(0.51)       \\   \hline
      1000  &    0.01     &      890        & 10000(0.34)       &     3240        &  10000(0.62)       \\   \hline
       50   &    0.001    &      688        &  5000(0.05)       &     $1898$      &   5000(0.29)       \\   \hline
      200   &    0.001    &     1980        &  5000(0.19)       &     6237        &   5000(0.51)       \\   \hline
      1000  &    0.001    &     2926        & 10000(0.34)       &     7859        &  10000(0.62)       \\
     \hline
   \end{tabular}
\end{center}
\end{table}

In the table we display the number of iterations of the Algorithm
\ref{alg1} (DSM) and Algorithm \ref{alg5} (PCM) for different values
of $m, M$ and of the accuracy $\epsilon$. For $m=50$ and $m=200$ the
maximum number of iterations that we allow is $5000$. For $m=1000$
we iterate  at most $10 000$ times.  When maximum number of
iterations is reached we also display between brackets the
corresponding accuracy.   We see that the duality gap is much better
with our Algorithm \ref{alg5} than with Algorithm \ref{alg1}.


\section{Conclusions}
A new decomposition
method in convex programming is  developed in
 this paper using the framework of dual decomposition.  Our method
 combines the computationally non-expensive gradient  method with the
 efficiency of structural optimization for solving
 separable convex programs. Although our  method
 resembles  proximal-based methods, it  differs both in the
 computational steps and in the choice  of the parameters. Contrary to most proximal-based methods that enforce  the
 next iterates to be close to the previous ones, our method  uses
 fixed centers  in the prox-terms  which leads to  more freedom in
 the next iterates. Another  advantage  of our proximal center method  is that it is
 fully automatic, i.e. the parameters  of the scheme are chosen
 optimally, which are crucial for justifying  its convergence
 properties.  We also presented
 efficiency estimate results of the new   method for
 general separable convex  problems and proved global convergence.
  The computations   on some test
 problems  confirm that the proposed method works  well in practice.


\section*{Acknowledgment}
\small{\textbf{Acknowledgment.}  We  acknowledge  financial support
by Research Council K.U. Leuven: GOA AMBioRICS, CoE EF/05/006,
OT/03/12, PhD/postdoc \& fellow grants; Flemish Government: FWO
PhD/postdoc grants, FWO projects
 G.0499.04, G.0211.05, G.0226.06, G.0302.07; Research communities (ICCoS, ANMMM, MLDM);
 AWI: BIL/05/43, IWT: PhD Grants; Belgian Federal Science Policy Office: IUAP
 DYSCO.}

\bibliography{kuleuven_bibtexfile}
\end{document}